\documentclass[fp,twocolumn]{jpsj3}

\usepackage{txfonts}
\usepackage{xurl}
\usepackage{graphicx}

\usepackage{listings}
\newcommand{\R}{\mathbb{R}}
\newcommand{\Z}{\mathbb{Z}}

\makeatletter
\newcommand\newblock{\hskip .11em\@plus.33em\@minus.07em}
\makeatother

\title{Persistent Homology Analysis for Materials Research and Persistent Homology Software: HomCloud}

\author{Ippei Obayashi$^{1,4}$\thanks{i.obayashi@okayama-u.ac.jp}, Takenobu Nakamura$^2$, and Yasuaki Hiraoka$^{3,4}$}
\inst{$^1$ Cyber-Physical Engineering Informatics Research Core (Cypher), Okayama University,
3-1-1 Tsushima-naka, Kita-ku, Okayama 700-8530 Japan \\
$^2$ National Institute of Advanced Industrial Science and Technology(AIST), 1-1-1, Umezono, Tsukuba, Ibaraki 305-8568, Japan \\
$^3$ Kyoto University Institute for Advanced Study. WPI-ASHBi. Kyoto University, Yoshida Ushinomiya-cho, Sakyo-ku, Kyoto 606-8501, Japan \\
$^4$ Center for Advanced Intelligence Project, RIKEN. 
1-4-1 Nihonbashi, Chuo-ku, Tokyo }

\abst{
This paper introduces persistent homology, which is a powerful tool to characterize the shape of data using the mathematical concept of topology. We explain the fundamental idea of persistent homology from scratch using some examples. We also review some applications of persistent homology to materials researches and software for persistent homology data analysis. HomCloud, one of persistent homology software, is especially featured in this paper.
}

\begin{document}
\maketitle

\section{Introduction}
This paper introduces the concept of persistent homology (PH), a data analysis technique based on the mathematical concept of topology. PH utilizes topological structures such as connected components, holes, rings, and voids to characterize the shape of data at multiple scales.

We often have the opportunity to analyze the geometric structure of materials science data. From atomic-scale data obtained by molecular dynamics simulations to larger-scale data obtained by various types of microscopy, materials science provides many spatial structures.
It is a typical problem to investigate how such structures are correlated with physical parameters and the physical properties of the materials.

We need to quantify the geometric structures to investigate the problems. PH enables us to summarize the shape of data quantitatively using mathematics.
From the viewpoint of materials informatics, PH provides descriptors of the shape of data.

The output of PH is a persistence diagram or a persistence barcode. 
The diagram and barcode have the same information, and only how to visualize is different. Therefore we will always use the persistence diagram in this paper.

This paper is organized as follows: Section \ref{subsec:recent} introduces the history and recent research trends in PH.
Section \ref{sec:ph} describes the theoretical foundations of PH, which will help the readers to understand how PH extracts the geometric features of data. Section \ref{sec:application} introduces the applications of PH to materials researches, which will help readers to understand what kind of data is suitable for PH.
Section \ref{sec:software} introduces the PH software, especially HomCloud, which is mainly developed by Ippei Obayashi, one of the authors.
Section \ref{sec:conclusion} summarizes this review paper and gives some concluding remarks.

\section{A Brief History of PH}\label{subsec:recent}

The concept of PH first appeared in the paper\cite{elz} by Edelsbrunner et al. in 2002. 
They studied topological persistence in a growing sequence of simplicial complexes in ${\bf R}^3$. 
Then, Zomorodian and Carlsson\cite{zc} introduced a mathematical framework of PH based on finitely-generated graded modules over polynomial rings ${\mathbf k}[t]$ with one variable $t$, where ${\mathbf k}$ is a field. This framework generalized the original concept of PH and enlarged the applicability of PH into wider classes of geometric problems. 
Furthermore, they clarified the so-called structure theorem of PH, i.e., under a certain finiteness condition, any PH can be uniquely decomposed into a direct sum of intervals (see Sect.~\ref{subsec:math}). This property provides a natural visualization of PH, called persistence diagrams, by plotting points with birth and death endpoints of intervals on the plane. The persistence diagrams show topological summaries of our input data in a multi-scale way. 

In order for PH to be applied to practical data analysis,  it should be clarified whether the persistence diagram is stable with respect to small perturbations. This property will be important for applying PH to practical problems since those data often contain noises or errors, and the output of PH should be stable to those small perturbations for obtaining essential data structure. The question of this type was first solved by Cohen-Steiner et al.\cite{ceh}, and they positively proved this property of persistence diagrams. Several generalizations of this result have also been reported later\cite{csggo,bl}. 

The structure and stability theorem explained here are the most important properties in applying PH to practical problems. The former provides a compact descriptor showing a topological multi-scale structure of data, and the latter guarantees their stability with respect to small perturbations. In fact, PH has been recently applied in various fields of science, including
sensor networks\cite{dg}, materials science\cite{hnhemn,stfrh}, 
biological evolution\cite{virus}, biomolecular structural analysis \cite{wei}, 
 brains\cite{brain}, cosmology\cite{cosmology},  etc. In those applications, a significant property of PH characterizing topological structures in a multi-scale way plays an important role in clarifying new insights. 
In this paper, we focus on the applications to materials science and show some of those successful examples. 

We also remark that these various applications give strong motivations for further mathematical studies of PH. For example, some methods of inverse problems and machine learnings of PH are actually developed by considering practical applications into materials science\cite{eh,linear-models}. 
For details of the recent development of PH-based machine learning and its applications, we refer to the survey paper\cite{ph_ml_survey}.

\section{Foundation of PH}\label{sec:ph}

This section describes how to characterize the shape of data quantitatively in the form of persistence diagrams.
PH is available for various kinds of data such as pointclouds and bitmaps, but we mainly consider pointclouds as input data in this paper.

\begin{figure*}[bt]
  \centering
  \includegraphics[width=0.9\hsize]{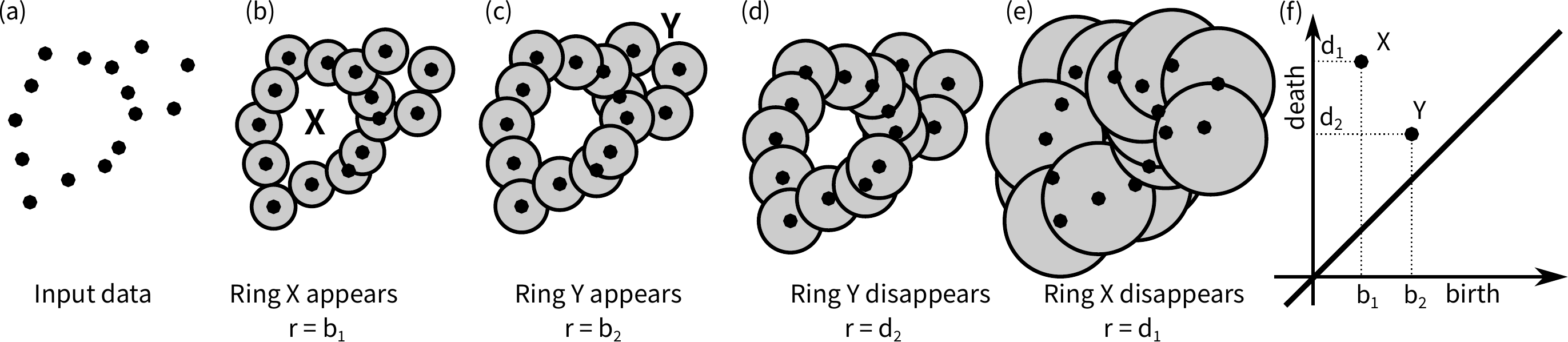}
  \caption{Filtration and PD (a) Input data (b) Pointcloud and discs with $r = b_1$ (c) $r=b_2$ (d) $r=d_2$ (e) $r=d_1$ (f) PD}
  \label{fig:ph}
\end{figure*}

\begin{figure}[tb]
  \centering
  \includegraphics[width=\hsize]{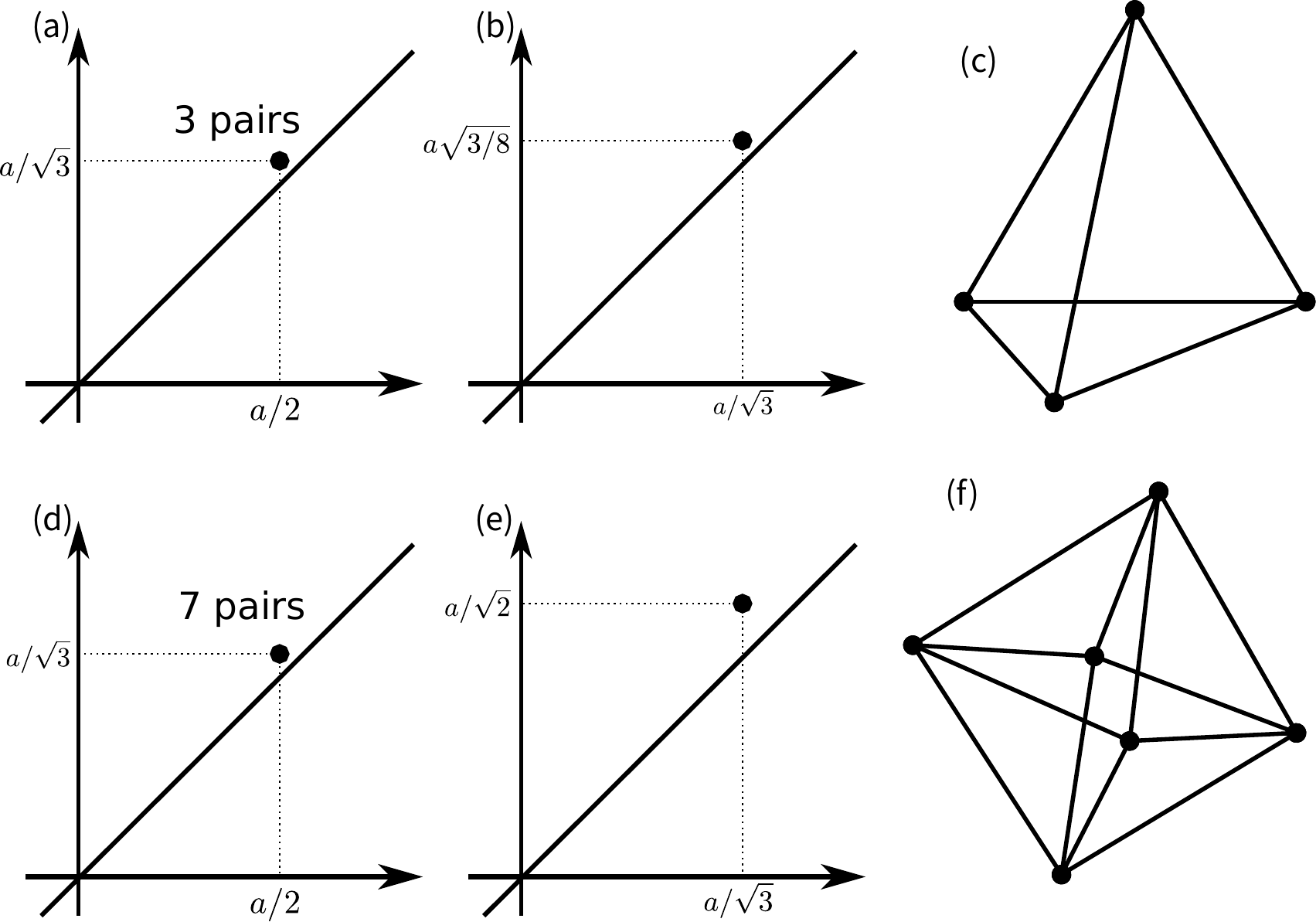}
  \caption{PDs for regular tetrahedral and octahedral points.
    (a) PD1 for tetrahedral points (b) PD2 for tetrahedral points (c) a regular tetrahedron
    (d) PD1 for octahedral points (e) PD2 for octahedral points (f) a regular octahedron
  }
  \label{fig:pd-tetrahedron}
\end{figure}

A pointcloud is a finite set of points. A typical pointcloud is an atomic configuration data.
We can analyze such data using PH.
The outline of PH is explained by the examples in Fig.~\ref{fig:ph},\ref{fig:pd-tetrahedron}.

A pointcloud Fig.~\ref{fig:ph}(a) has no ring or hole, but it seems to have two rings.
To construct topological structures on the data, we put discs on all points whose radii are the same, and holes appear as in Fig.~\ref{fig:ph}(b-e). We can count the number of holes in these figures.

The problem here is to determine the radius of discs. The number of holes changes when the radius changes.
The fundamental idea of PH is to consider the changing process instead of fixing the radius.
When the radius $r$ gradually becomes larger from zero, holes appear and disappear as in Fig.~\ref{fig:ph}(b-e).
The increasing process of the shapes is called a \emph{filtration}, and the theory of PH makes pairs of appearance and disappearance of holes.
In Fig.~\ref{fig:ph}, $(b_1, d_1)$ and $(b_2, d_2)$ are paired.
The radius of appearance is called a \emph{birth time}, the radius of disappearance is called a \emph{death time}, and the pair of birth and death times is called a \emph{birth-death pair}.
The set of birth-death pairs with multiplicity is called a \emph{persistence diagram} (PD).
The PD is often visualized by a scatter plot or 2D histogram (Fig.~\ref{fig:ph}(f)).
We remark that the birth times are death times are sometimes squared according to conventions in computational geometry and topological data analysis.

PH is applicable to any dimensional data. For 3D data, we use spheres instead of discs. 
In homology theory, 3D geometric objects have three types of homology information, and each type is characterized by dimension.
0D homology has information about connectivity, 1D homology has information about holes or rings, and 2D homology has information about cavities or bubbles. Corresponding 0D, 1D, and 2D PDs are available. Since we consider rings in Fig.~\ref{fig:ph}, in Fig.~\ref{fig:ph}(f) is a 1D PD.

One advantage of PH is its mathematical background. 
The structure theorem of PH gives the algorithm to compute the diagram, and the theorem also ensures the uniqueness of the diagram.
It means that the same input gives the same diagram, unlike the Monte-Carlo method or stochastic gradient descent.
The stability theorem of PH ensures that the small change of the input data causes only small changes of output PD.
These theorems play an important role in reliable data analysis.

To understand PDs, we examine some PDs for typical pointclouds. 
Figures~\ref{fig:pd-tetrahedron}(a) and (b) show the PDs of regular tetrahedral points (Fig. \ref{fig:pd-tetrahedron}(c)) where $a$ is the distance between two points.
PD1 for regular tetrahedral points has three birth-death pairs, and all of them are $(a/2, a/\sqrt{3})$.
PD2 has one birth-death pair $(a/\sqrt{3}, a\sqrt{3/8})$.
Three pairs in PD1 correspond to triangles of the tetrahedron, and one pair in PD2 corresponds to the void at the center of the points.
$a/2$ is half the length of the edge, $a/\sqrt{3}$ is the circumradius of the triangle, and $a\sqrt{3/8}$ is the
circumradius of the tetrahedron.

Probably readers will wonder why PD1 has three pairs, not four. The reason comes from the mathematical theory of homology. We will discuss the background in Sect.~\ref{subsec:math}.

We can apply the same idea to the regular octahedral points.
Figures~\ref{fig:pd-tetrahedron}(d) and (e) show the PDs for the pointcloud.
PD1 for the regular octahedral points has seven pairs at $(a/2, a/\sqrt{3})$.
These pairs correspond to the triangles of the octahedron. The PD1 does not have pairs corresponding to the eighth triangle and
the square for the same reason as the tetrahedron. PD2 has one birth-death pair at $(a/\sqrt{3}, a/\sqrt{2})$ corresponding to
the void at the center of the octahedron.

The PDs of hcp crystalline structure are in some sense the combination of PDs for tetrahedral points and octahedral points.
PD1 for the hcp has only birth-death pairs at $(a/2, a/\sqrt{3})$ corresponding to the triangles, and
PD2 has birth-death pairs at $(a/\sqrt{3}, a\sqrt{3/8})$ and $(a/\sqrt{3}, a/\sqrt{2})$ corresponding to
the tetrahedral sites and octahedral sites in the hcp structure. Indeed fcc crystalline structure also has the tetrahedral sites and octahedral sites, and the PDs for the fcc points are completely the same as hcp.

\subsection{Geometric representations}

The concept of PH is applicable to various data.
We now describe some geometric representations to understand available data for PH.

One important geometric representation is called a \emph{simplicial complex}.
A simplicial complex is composed of points, line segments, triangles, tetrahedrons, and higher dimensional counterparts.

Formally saying, an $n$-dimensional simplicial complex is a finite set of $k$-simplices for $k=0,\ldots,n$, which are represented by $k+1$ vertices.
A 0-simplex is a point, 1-simplex is a line segment, 2-simplex is a triangle, and 3-simplex is a tetrahedron.

For an $n$-simplex $\sigma$ and $k < n$, a $k$-simplex included in $\sigma$ is called a face of $\sigma$ if the all vertices of the $k$-simplex is also the vertices of $\sigma$.

A finite set of simplices $X$ is called a \emph{simplicial complex} if $X$ satisfies the following two conditions:
\begin{enumerate}
\item If $\sigma \in X$, any face of $\sigma$ is contained in $X$
\item The intersection of two simplices in $X$ is their common face
\end{enumerate}
Condition (2) means that two simplices are glued together by a lower-dimensional simplex.

A purely combinatorial description of a simplicial complex is called an abstract simplicial complex.
For a finite set $V$, a family of subsets of $V$, $\Sigma$, is called \emph{abstract simplicial complex} if the following two conditions hold:
\begin{itemize}
\item For any $v \in V$, $\{v\} \in  \Sigma$
\item If $\sigma \in \Sigma$ and $\tau \subset \sigma$, $\tau \in \Sigma$
\end{itemize}
Of course, We can regard a simplicial complex as an abstract simplicial complex.
An abstract simplicial complex is helpful to represent a geometric object on a computer.

One important simplicial complex is an \emph{alpha complex}\cite{alpha}.
We can construct an alpha complex with radius parameter $r$ from pointcloud as a subset of Delaunay triangulation using the Voronoi diagram.
An alpha complex with parameter $r$ must have the same topological information as the union discs model shown in Fig~\ref{fig:ph}.
Therefore we usually use an alpha complex to represent the union discs model such as Fig~\ref{fig:ph}(b-e).

\begin{figure}[tbp]
  \centering
  \includegraphics[width=\hsize]{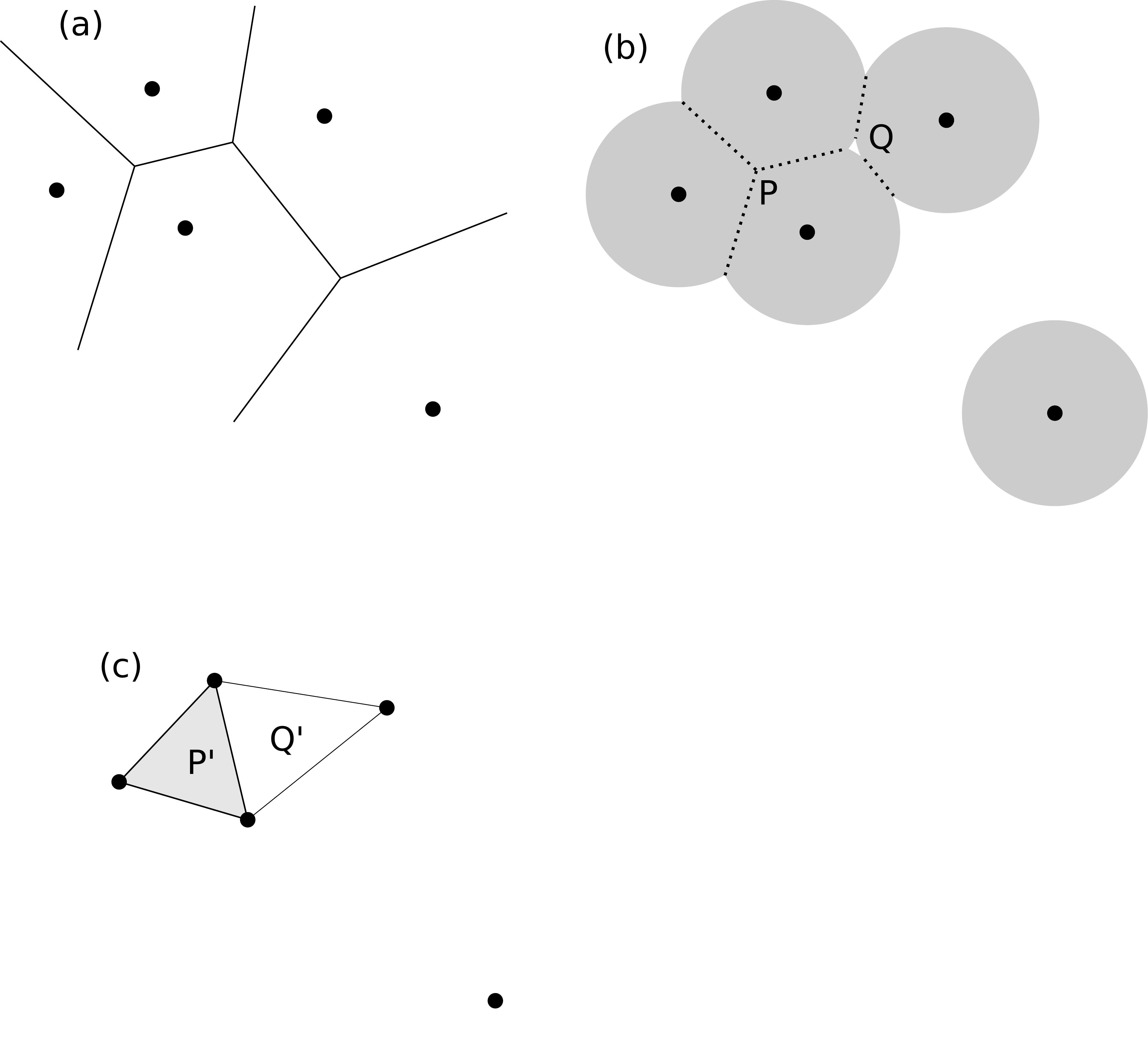}
  \caption{Construction of an alpha complex (a) Pointcloud and Voronoi diagram (b) Overlay image of discs and Voronoi diagram (c) Alpha complex}
  \label{fig:alpha}
\end{figure}

We explain the construction of an alpha complex using the example shown in Fig.~\ref{fig:alpha}.
Figure~\ref{fig:alpha}(a) shows a pointcloud and its Voronoi diagram, (b) shows the overlay image of (a) and discs, and (c) shows the alpha complex of the pointcloud.
The triple intersection $P$ in (b) corresponds to the triangle $P'$ in (c), and the double intersections in (b) (shown in dotted lines) correspond to edges in (c).
Since there exists no triple intersection at $Q$ in (b), the alpha complex has no triangle at $Q'$ in (c).
As shown in Fig.~\ref{fig:alpha}, (b) and (c) have the same topological information.
For example, both (b) and (c) have one hole ($Q$ and $Q'$) and two connected components.
Nerve theorem from abstract topology ensures that the union of discs and the alpha complex gives the same topological information.
A filtration by alpha complexes is called an \emph{alpha filtration}.

We can extend the idea of the alpha complex to 3D or higher dimensional spaces by considering quadruple intersections and more.

Another important simplicial complex for data analysis is a clique complex and Vietoris-Rips complex.
We can define a clique complex from an undirected graph.
For a graph $G = (V, E)$, a set of vertices $C = \{v_0, \ldots, v_k\}$ is called a \emph{clique} if any pair $(v_i, v_j)$ is adjacent, that is, the graph $G$ has an edge between $v_i$ and $v_j$.
The set of all cliques is called a \emph{clique complex}. 
We easily show that the clique complex satisfies the conditions of an abstract simplicial complex, and we write it as $X(G)$.
By using clique complexes, we can construct filtration from a weighted graph. For a weight $w: E \to \R$ and $a \in \R$, $G_a = (V_a, E_a)$ is a subgraph of $G$, where:
\begin{equation}
  \begin{aligned}
    V_a &= V, \\
    E_a &= \{ e \in E \mid w(e) \leq a \}.
  \end{aligned}
\end{equation}
$X(G_a)$ gradually increases when $a$ increases, so $\{X(G_a)\}_{a \in \R}$ gives a filtration of clique complexes.

By considering a finite number of points $V$ on a metric space, we can construct a filtration called a \emph{Vietoris-Rips filtration} using the complete graph on $V$ and the weight function $w(\{v_i, v_J\}) = d(v_i, v_j)$, where $d(v_i, v_j)$ is the distance between two points.

Vietoris-Rips filtrations are sometimes used for PH data analysis.
In the study of molecular phylogenetics by Chen et al.\cite{virus}, they use Vietoris-Rips filtrations using a genomic distance between genomic sequences.

One advantage of Vietoris-Rips filtrations is that they can be used as long as pairwise distances are available.
At the same time, the disadvantages are high computation cost and lack of geometric correspondence as in alpha filtrations.

We also use cubical complexes to analyze pixel and voxel data. 
A cubical filtration consists of points, line segments, squares, cubes, and higher-dimensional counterparts.
We can construct a filtration of cubical complexes by considering a level function on pixels.
The concept can be directly applicable to gray-scale bitmap data.
This type of filtration is called ``level-set filtration''.
Some studies\cite{iron-ore,porous,capillary-fluid,flow-estimation} apply a cubical filtration to 2D or 3D binary bitmap data using distance transform. 

\subsection{Mathematical foundation}\label{subsec:math}

PH is based on the theory of homology, a part of topology theory. 
Since this paper is intended for readers who are interested in the application of PH, we do not attempt to fully formulate the mathematics of PH.
Instead, we will illustrate the mathematical ideas using some examples.

\begin{figure}[tbp]
  \centering
  \includegraphics[width=0.7\hsize]{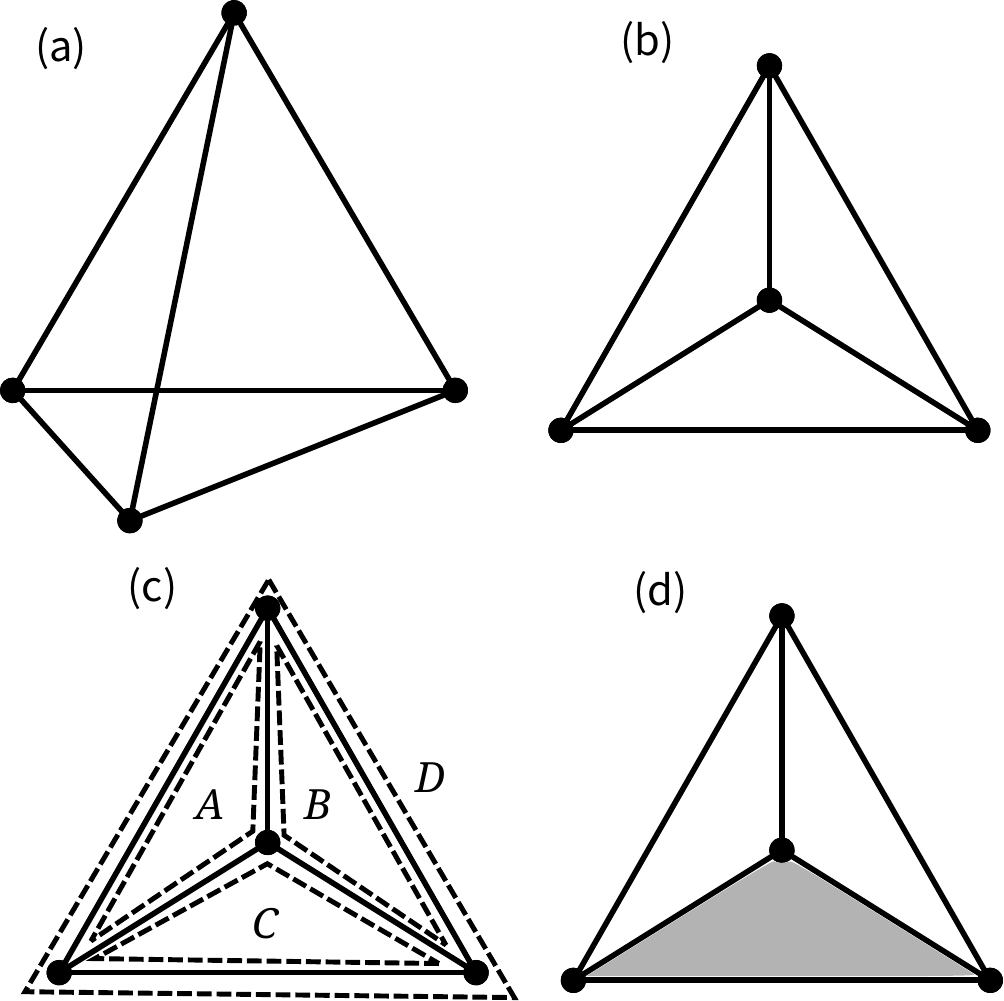}
  \caption{(a) Tetrahedron (b) Tetrahedron looking from above
    (c) Four rings in the tetrahedron
    (d) A ring $C$ in (c) is filled with a triangle}
  \label{fig:tetrahedron}
\end{figure}

First, we consider the problem of counting the number of rings in a tetrahedral skeleton.
The tetrahedron (Fig.~\ref{fig:tetrahedron}(a)) seems to have four rings, but it looks like three rings when we see the tetrahedron from above (Fig.~\ref{fig:tetrahedron}(b)).
This is because the outer ring $D$ in Fig.~\ref{fig:tetrahedron}(c) looks like the combination of the three rings $A$, $B$, and $C$ in some sense.

The theory of homology justifies intuition using linear algebra.
In homology theory, each ring is regarded as a vector, and we can justify the equality $A + B + C = D$ by assuming the rule that adding the same edge twice equals zero.
An easy way to justify the rule is to consider $\Z/2\Z$-vector space.
The equality means that the four rings $A, B, C, D$ are not linearly independent.
We can consider the linear space of all rings, and we can count the number of linearly independent rings by computing the dimension of the linear space.
In Fig~\ref{fig:tetrahedron}(a)(b), the tetrahedron skeleton has three linearly independent rings.
This is because PD1 for regular tetrahedral points has three birth-death pairs at $(a/2, a/\sqrt{3})$, not four.
The same is true for the regular octahedral points.

Next, we consider the counting of holes in Fig.~\ref{fig:tetrahedron}(d).
The figure seems to have two holes since the hole $C$ in Fig.~\ref{fig:tetrahedron}(c) is filled with a triangle and the rest are $A$ and $B$.
In other words, we can count the number of holes by $z - b$, where $z$ is the number of linearly independent rings, and $b$ is the number of linearly independent rings filled by triangles.
We can realize the idea by canceling out the ring $C$ in an algebraic way called the quotient.
By quotient, we can compute the linear space of all unremoved rings.
In fact, $z - b$ is the dimension of the quotient linear space.
Of course, the same idea is available to count the number of cavities. 

\begin{figure*}[bt]
  \centering
  \includegraphics[width=0.6\hsize]{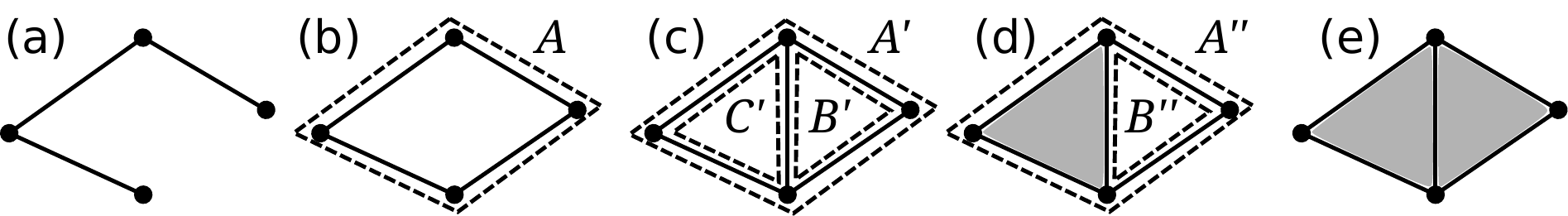}
  \caption{Increase sequence of simplicial complexes}
  \label{fig:phmath}
\end{figure*}

Finally, we explain the idea of PH using the example in Fig.~\ref{fig:phmath}.
In this figure, the numbers of holes are $0, 1, 2, 1, 0$ from left to right.
The sets of linearly independent rings without filling with triangles are $\{\}, \{A\}, \{B', C'\}, \{B''\}, \{\}.$
To understand the relationship between these sets, we consider the change of basis.
By changing the basis from $\{B', C'\}$ to $\{B' + C', C'\}$, we can consider the following relationship called \emph{interval decomposition}, where $*$ means the vanishment of the hole:
\begin{equation}\label{eq:decomposition}
  \begin{array}{ccccccccc}
    \{\}&\to&\{A\}&\to&\{B' + C', C'\}&\to&\{B''\}&\to&\{\}  \\ \hline
        && A & \to &A' = B' + C' & \to & A'' = B'' &\to &* \\
        &&&& C' &\to &*
  \end{array}
\end{equation}
In the above relationship, $ A'' = B''$ is justified due to canceling out rule.
We can say that the hole $A$ appears at Fig.~\ref{fig:phmath}(b) disappearance at Fig.~\ref{fig:phmath}(e), and another hole $C'$ appears at Fig.~\ref{fig:phmath}(c) and disappears at Fig.~\ref{fig:phmath}(d).
We can make pairs of appearance and disappearance of the holes mathematically in that way.
The theory of PH ensures that we can always find such proper bases.
This is the fundamental idea of PH.

\section{Applications to Materials Research}\label{sec:application}

PH, especially for point clouds, is currently being used as a tool to extract structures in various materials.
In Ref. \cite{nakamura2015persistent}, two of the authors already give some case studies to show how PH works in representative structures found in material systems.
Since then, many interesting applications have been investigated, including granular materials, network former, and polymer.
There are several properties of these systems that explain the effectiveness of PH feature extraction.
They are not found in the introduction of the previous case studies.
Now, we will focus on applications of the point clouds and present both the validity of using PDs as a descriptor for these materials and the validity of geometric interpretations for them.
We then will present some notable recent results for each of these materials.

\subsection{PD as a descriptor of disordered system}

PH has the potential to become a universal tool for representing complex structures in material science.
One of the central issues in materials science is to express the difference in physical properties by the difference in structure.
This approach has been very effective for regularly arranged structures, such as crystals, and extremely disordered structures, such as gases.
Apart from these two extremes, however, it is difficult to take such an approach because of the difficulty in quantitatively evaluating complex structures.
For example, there can be structures partially ordered in some regions while disordered in other regions, or structures in an intermediate stage of ordering. 
To quantify them, the method must be able to express the diversity of the structure.
Since PD can be regarded as a kind of distribution function, it can represent those features.

When extracting a characteristic shape from a disordered configuration of particles, we may face the problem of the arbitrariness of the threshold, as mentioned in Sect. 3 as the problem of determining the radius of the disk. 
For example, in order to introduce hydrogen bonds, metal complexes, polyhedra of metal clusters, etc., it is necessary to determine an appropriate threshold.
Usually, thresholds are set based on physical and chemical properties, but PH provides an alternative solution to this problem.
PH evaluates the degree of robustness of the thresholds by sweeping them.
Holes with a small $d-b$, that is, holes that disappear soon after they appear as the radius parameter increases, are not detected unless the threshold is finely tuned, and vice versa.
In other words, holes far from the diagonal $b=d$ in the PD can be extracted even if some numerical noise is added.
Therefore, we can assume that such holes have important information for the disordered configuration.

\subsection{Comparison with conventional methods}
To gain a better understanding of the uniqueness of PH, it will be helpful to compare it with conventional methods used in material science such as Voronoi analysis\cite{cheng2011atomic}, ring statistics\cite{zallen2008physics}, bond orientational order parameter\cite{steinhardt1983bond}, and radial distribution function\cite{hansen1990theory}.
The Voronoi analysis and the ring statistics provide discrete variable indicators.
Due to the discreteness, two configurations close to each other may be identified as different indices.
In PH, the stability theorem allows us to avoid this problem.

The bond orientational order parameter is a continuous variable indicator.
Together with the two mentioned above, this indicator expresses the degree of disorder in terms of similarity to the crystalline structure.
When the structure of interest is relatively close to the crystal, these quantities work for interpreting the structure.
However, they are not suitable for representing disordered shapes far from the crystalline structure.
Whereas, PH provides interpretations that are not based on similarity to crystals.

The radial distribution function is a quantity that does not rely on the similarity of crystal and can be applicable to a disordered structure.
In fact, it can completely express the disordered structure of the simple liquid.
However, it is not suitable for expressing complex shapes composed of many particles, which can be expressed by PH.

As can be seen from these comparisons, 
PH is suitable for describing shapes consisting of a large number of particles embedded in a structure that is far from the crystal structure.
Then, the degree of the disorder can be expressed as a distribution function.
A material system that has the benefit of applying PH has these characteristics.

\subsection{An object with disordered structure}
We can get an overview of the systems in which PH is currently successfully applied by analyzing the reasons for the realization of disordered structure from the viewpoint of material science.
There are three possible reasons for the disordered structure:
1. because of the property of the microscopic constituents of the material,
2. because of the property of the process of preparing the material,
3. because the material is composed of macroscopic components with heterogeneity and friction.
Here we will describe them.

Soft matter, such as polymers, often has disordered structures on the atomic and molecular scales.
In molecular systems that constitute soft materials, there are intermolecular potential and intramolecular potential such as a torsional potential. The magnitude of these energies is of the same order as the thermal energy at room temperature.
As a result, both single molecules and molecular assemblies are allowed to have various structures in soft matter and take quite disordered structures.

The solid-state with a disordered structure, called glass or amorphous, is obtained by quenching the material while maintaining the disorder that appeared in the liquid state.
Many multi-component systems reach the glassy state by this procedure.
In addition to the soft matter, alloys and ceramics can also achieve a glassy disordered structure.

Granular consisting of macroscopic particles also often have disordered structures.
Unlike atoms or molecules, the granular particles have diversity in size and interactions.
Therefore, when they are densely packed, it is not necessarily arranged regularly like crystals.
Even if particles have a uniform property in the simulation, they may not always form crystals and solidify in a disordered structure due to friction.
Therefore, in both experiments and simulations, granular systems often show a disordered structure.

For these reasons, these material systems often have structures that are difficult to quantify, and PH is required to be quantified there.

\subsection{Material systems with geometric interpretation}
We can interpret the shape extracted by PH  if the physical process of creating the point cloud is governed by geometry. 
In other words, the validity of interpreting the point cloud as a complex consisting of an expanding sphere accompanying each point is justified on the basis of materials science.
For a mono-disperse system, the construction is always justified, but for a multi-component system, it is not always.
It is non-trivial that each particle has an intrinsic quantity corresponding to its radius. 
Only if this property is present can we define a physically interpretable contact between a particle and another adjacent particle.
In a system with $N$ kinds of particles, there are $N(N+1)/2$ kinds of inter-particle distances.
Since there must be $N$ radii of particles, the degree of indefiniteness is $N(N-1)/2$.
The radius is well defined only when this indefiniteness is not present.

Since the configuration of the system is determined by the inter-particle potential between the particles, let us consider the conditions under which the potential has the distance corresponding to the radius.
For example, consider a binary system ($N=2$) with a Lennard-Jones potential $U_{AB}(r)=\epsilon_{AB}[(\sigma_{AB}/r)^{12}-(\sigma_{AB}/r)^6]$, where $A$ and $B$ represent the type of particle.
The length of the interaction between component A and component B is determined by $\sigma_{AB}$.
As mentioned above, there are three length parameters, $\sigma_{AA},\sigma_{AB},\sigma_{BB}$, corresponding to the type of distance between the particles.
If $\sigma_{AB}$ is expressed using the arithmetic mean of $\sigma_{AA},\sigma_{BB}$, we can introduce the radius as $\sigma_A=\sigma_{AA}$ or $\sigma_B=\sigma_{BB}$.
Here, the condition of the arithmetic-mean has eliminated the indefiniteness for the radius.
The same applies to the system with $N$ kinds of particles.
Especially in the case of LJ, the equation given by this arithmetic mean is called the Lorentz-Bertelot combining rule\cite{delhommelle2001inadequacy}.
Thus, the combining rule guarantees the validity of considering the radius of the particle and makes the virtual contact in alpha shape physically meaningful.

Many polymer systems currently analyzed might be considered to satisfy the combining rule.
The interaction between granular particles always satisfies it whether in an experimental or simulated system because the particle is macroscopic.
The details will be discussed later respectively.
In the case of network former, radius and contact can be reasonably introduced from another justification.
This will also be discussed later.

\subsection{Achievements in material systems}
In the following, we will introduce some notable recent results of PH in the material system.
Most of the analyses are point clouds generated by particle simulations, and some data are generated by processing data from the inverse Monte Carlo method based on the experimental measurement.

\subsubsection{Granular}
Granular systems have been studied extensively since the very early days of PH's application to material systems.
Unlike the atomic and molecular scales, granular particles are macroscopic objects, and the radius and shape are well defined without ambiguity.
In fact, the interactions between granular particles are usually represented by the contact force such as repulsive Hertzian or harmonic.
These potentials satisfy the combining rule, which guarantees that PH works effectively.
Concerning extraction of structure, local structures close to crystalline structure and deformation with geometric constraints have been studied \cite{stfrh}. 

In addition to the extraction of structure, PH also helps us to quantify the force-chain network.
Due to the contact force, the mechanical properties at the particle scale are characterized by a force chain network.
This is a unique property to granular systems.
The bulk mechanical properties are then found to be related to the Betti number of the force chain network.
Reference \cite{kondic2012topology} is one of the very early studies of PH in material systems, which found that friction and polydispersity are expressed as differences in the force chain network and that the force chain network changes qualitatively at the Jamming transition point.

PH is also useful for expressing the process dependency of the force chain network.
The static mechanical properties of granular depend on the process.
This subject is one of the central issues of research of granular systems, it is often addressed in TDA.
In Ref.\cite{ardanza2014topological}, it was found that PD1 can represent the history dependence of the tapping operation of the granular system even for the system with the same density.
This is in contrast to liquids, where density is a state quantity that represents the physical properties of the system.
In addition to the static properties, the transient properties are also characterized by PH.
Changes in the force chain network associated with the impact have been quantified for both the granular system  \cite{lim2018topology}  and the suspension system \cite{hayakawa2021impact}.

\subsubsection{Network former}
Network former, such as silica, is another material where PH is actively applied.
In the case of network former, the combining rule may not hold for the interaction parameters.\cite{coslovich2009dynamics}
However, there is another reason that justifies interpreting the extracted shapes. 
Since the closest particle type to each particle is limited, a radius can be introduced for each particle species.
In the realized particle configuration, indefiniteness to determine radius does not practically appear because the type of particle with the closest coordination to each particle is restricted.

As an example, let us consider the silica (${\rm SiO}_2$) system. 
In the case of silica at low temperatures, silicon(Si) is always at the center of a tetrahedral structure consisting of oxygen(O), and no two silicons are directly next to each other.
Therefore, the types of inter-particle distances are limited to Si-O and O-O, and Si-Si does not appear. 
The radius of the oxygen can be determined from the nearest neighbor distance of O-O, and the radius of Si can be determined from the nearest neighbor distance of Si-O. 
Consequently, we can interpret the shape realized by silica through the union of balls that appeared in PH.
If there are additional elements, there is no guarantee that the interpretation will work.
However, we may ensure the validity of the interpretation by considering additives as secondary effects added to the backbone structure.

PH was first applied to network former by Hiraoka et al.\cite{hnhemn}
It has been pointed out that the PD1 is associated with a structure so-called medium-range order, and the correspondence with the length scale of the first sharp diffraction peak has been discussed.
By focusing on the property that the curves appearing in PD1 correspond to geometric constraints, they found a hidden structure in the disordered configuration.
A detailed interpretation of the geometric constraints has now been illustrated using a simplified model \cite{ormrod2021persistent}.
Moreover, it was found that PD has the ability to estimate the glass transition temperature using the support vector machine\cite{kusano2017kernel}.
PH has also been applied to experimentally obtained configuration data for various network formers\cite{Onodera2020,Koyama2020}.
Amorphous ice is another example of a network former, where the medium-range order of the hydrogen-bonding network is quantified by PH.\cite{hong2019medium}

\subsubsection{Polymer}
Polymers are another material where PH is currently being actively applied.
As mentioned earlier, soft matter systems are very compatible with PH analysis because of the need to extract characteristic shapes from the disordered configuration.
In many cases, coarse-grained models are used in molecular simulations of polymers.
In general, there is no guarantee that the combining rule holds for the interaction potential of coarse-grained models.
For example, the combining rule breaks down when charged particles or heterogeneous particles are treated as a single particle in coarse-graining.
However, since the system currently under analysis does not meet these conditions, such contribution is expected to be small.
Therefore, the PH interpretation can be applied to data obtained by coarse-grained molecular dynamics.

In the case of polymer systems, the correlation with material properties is actively discussed rather than the extraction of the structure itself, compared to the two examples mentioned above.
For example, correlations with dielectric constant \cite{Shimizu2021}, shear response \cite{yoshimoto2021molecular}, and crazing process \cite{ichinomiya2017persistent} were discovered.
As a study that focuses on the structure itself, a method to quantify the threading of ring polymers has been proposed \cite{PhysRevResearch.2.033529}.

\subsubsection{Other systems}

In addition to these representative systems, PH has also been applied to other material systems.
In ${\rm Pd}_{40}{\rm Ni}_{40}{\rm P}_{20}$ bulk metallic glasses, the medium-range ordered structure was extracted by applying PH to the configuration data of both total and each component \cite{hosokawa2019partial}.
For the toy model of amorphous solids by simulation, the effectiveness of PH in describing the structural changes during plastic deformation was evaluated based on machine learning \cite{rocks2021learning}, the yielding was associated with a decrease in the number of robust holes introduced by PH.\cite{shirai2019microscopic}.
In the experiment, PH was used to express the local structure of two-dimensional binary colloidal configuration confined at the gas-liquid interface.\cite{doi:10.1021/acs.langmuir.8b01411}.
In addition to point clouds, PH has been applied to the visualization of energy landscape \cite{carr2016energy} and characterization of spatial patterns such as phase separation structures of magnetic materials\cite{Kotsugi} and polymers \cite{published_papers/31169608}. 

\subsection{Perspectives for material science}

We have presented some recent applications of PH, mainly to granular materials, network formers, and polymers.
It is expected that PH will be applied to other material systems in the future.
For example, soft matter other than polymers, such as liquid crystals, emulsions, vesicles, and colloidal gels, will be promising applications of PH because of the need to extract characteristic structures embedded in disordered configurations.
In addition to PH for point clouds, it is also promising to use PH for the quantification of systems with bubbles, fillers, and porous shapes because their shapes are dominated by holes.

\section{Software}\label{sec:software}

For the applications of PH, software is important.
The development of algorithms and software has progressed in parallel with the development of theories.
The paper by Edelsbrunner et al.\cite{elz} showed an algorithm to compute a PD, and the algorithm
has been refined theoretically and practically by subsequent researches\cite{zc,twist,cohomology1,cohomology2,eirene,multifield}.
Parallel and distributed algorithms has been also studied\cite{chunks,dipha}.

Many researchers on PH have developed various data analysis software using PH,
including 
Javaplex (\url{http://appliedtopology.github.io/javaplex/}),
Perseus\cite{perseus} (\url{http://people.maths.ox.ac.uk/nanda/perseus/}),
PHAT\cite{phat} (\url{https://bitbucket.org/phat-code/phat/}),
Dipha\cite{dipha} (\url{https://github.com/DIPHA/dipha}),
Ripser\cite{ripser} (\url{https://github.com/Ripser/ripser}),
Gudhi (\url{https://gudhi.inria.fr/}),
Dyonisys (\url{https://mrzv.org/software/dionysus/}, \url{https://mrzv.org/software/dionysus2/}),
R-TDA (\url{https://cran.r-project.org/package=TDA}),
EIRENE\cite{eirene} (\url{http://gregoryhenselman.org/eirene/}),
CubicalRipser (\url{https://github.com/CubicalRipser}),
RIVET (\url{https://github.com/rivetTDA/rivet}),
giotto-tda\cite{giotto-tda} (\url{https://github.com/giotto-ai/giotto-tda}),
jHoles\cite{jholes}, and
HomCloud (\url{https://homcloud.dev}).

The developers have their interests and analysis targets.
For example, Ripser focuses on the efficient algorithm to compute PDs for Vietoris-Rips filtrations.
A benchmark\cite{roadmap} showed that Ripser is one of the fastest software compared to competitors.
PHAT and Dipha are developed by the same people as Ripser, and they also have good performance.

Gudhi is more interested in computational geometry.
Gudhi collaborates with CGAL (Computational Geometry Algorithms Library, \url{https://www.cgal.org/}), a famous computer geometry library, and works with CGAL on the development.
Gudhi has various representations of geometric objects.

R-TDA provides an interface from the R language. In fact, R-TDA uses Gudhi, PHAT, and Dyonisys as a backend. R-TDA is a bridge between R and the TDA world.

Each software has its advantages and features.

\subsection{HomCloud}\label{subsec:homcloud}

In this subsection, we introduce our software, HomCloud. Ippei Obayashi, one of the authors of this paper, mainly develops HomCloud.

HomCloud is free software, and you can download HomCloud from \url{https://homcloud.dev}. You can freely use, copy, modify, and
redistribute the software
under GPL (\url{https://www.gnu.org/licenses/gpl-3.0.html}).

HomCloud focuses on applications, mainly to materials science.
We use HomCloud to analyze atomic configuration given by molecular dynamical
simulations and reverse Monte-Carlo and pixel data and voxel data given by
an electron microscope. Of course, other data than materials science is also available.
HomCloud has been already used in various scientific researches, including materials science\cite{iron-ore,Hirata2020,YoheiONODERA201919143,Onodera2020,Shimizu2021,Koyama2020,Onodera2019a,PhysRevB.99.045153,doi:10.1098/rspa.2020.0170,ANDO2021142112,Hong_2019,PhysRevResearch.2.033529,minamitani2021topological,yoshimoto2021molecular}, geology\cite{flow-estimation,SUZUKI2020104550}, structural biology\cite{ICHINOMIYA20202926}, and medical image analysis\cite{Koseki2020,Oyama2019}.

HomCloud has useful functionalities such as visualization, inverse analysis,
machine learning. Especially, inverse analysis of HomCloud, which detects
a ring or a cavity corresponding to each birth-death pair, is the most
advanced among other software.

HomCloud has two types of interface, command-line interface and
python interface. Since Python has a rich scientific computing ecosystem,
you can combine the output of HomCloud with the ecosystem
using the Python interface.

HomCloud is available on Windows, Linux, and macOS, including Apple Silicon Mac.
You can also use HomCloud on Google Colaboratory, which allows you to
write and execute Python code in the web browser. You can try HomCloud
on Colaboratory without installing HomCloud on your machine.

To improve the software quality, the developers of HomCloud do several practices.
One practice is dogfooding; that is, the developers use HomCloud
for daily data analysis.
Close communication between developers and users is effective for software
improvement, and dogfooding is one extreme way.
Another practice is continuous integration.
After a code is uploaded to the code repository, HomCloud is automatically built and tested.
Continuous integration is essential to keep the quality and portability of the software.

HomCloud internally uses various third-party components to reduce the development cost.
For example, HomCloud uses PHAT, Ripser, and Dipha to compute PDs.
These are known for their good performance; therefore we use them.
CGAL is also used to compute alpha filtrations.
Python's standard scientific computing libraries, such as
NumPy, SciPy, and Matplotlib are also used.

\subsubsection{Installing HomCloud}

Since HomCloud is written in Python, you need to install Python before
installing HomCloud. After installing Python, 
you can easily install HomCloud using pip command, which is
the de facto standard package management system for Python.
We recommend pip to install HomCloud on Linux, Intel macOS, and Windows.
Another installation option is conda (\url{https://conda.io}).
Conda is the open-source
package-management and environment management system mainly developed by
Anaconda Inc. (\url{https://www.anaconda.com}). Conda is a famous tool among
data science people. HomCloud prepares conda packages for Windows, Linux,
and Apple Silicon Mac.  We recommend conda for Apple Silicon Mac since conda-forge has the most
extensive package for Apple Silicon Mac. 
The installation manual is available at
\url{https://homcloud.dev/install-guide/index.en.html}.

\subsubsection{Basic input/output}
HomCloud accepts the following data.
\begin{itemize}
\item 2D/3D pointcloud (alpha filtrations are used)
\item Any dimensional bitmap data, both binary data, and grayscale data
  (cubical filtrations are used)
\item Distance matrix (Vietoris-Rips filtrations are used)
\item Abstract simplicial complex with weights 
\end{itemize}
We can assign an initial radius to each particle for a 3D pointcloud
to reflect the size of the particle using
weighted alpha shape\cite{weighted-alpha}.
This functionality is useful to reflect the physical radii such as
ionic radii.
We can also use periodic boundary conditions for 3D pointcloud.

The output of HomCloud is a PD.
You can plot a histogram by HomCloud.
The backend of plotting is Matplotlib, so you can create fancy figures by utilizing Matplotlib.
You can also output the list of birth and death times.

\subsubsection{Inverse analysis}

\begin{figure}[btp]
  \centering
  \includegraphics[width=0.6\hsize]{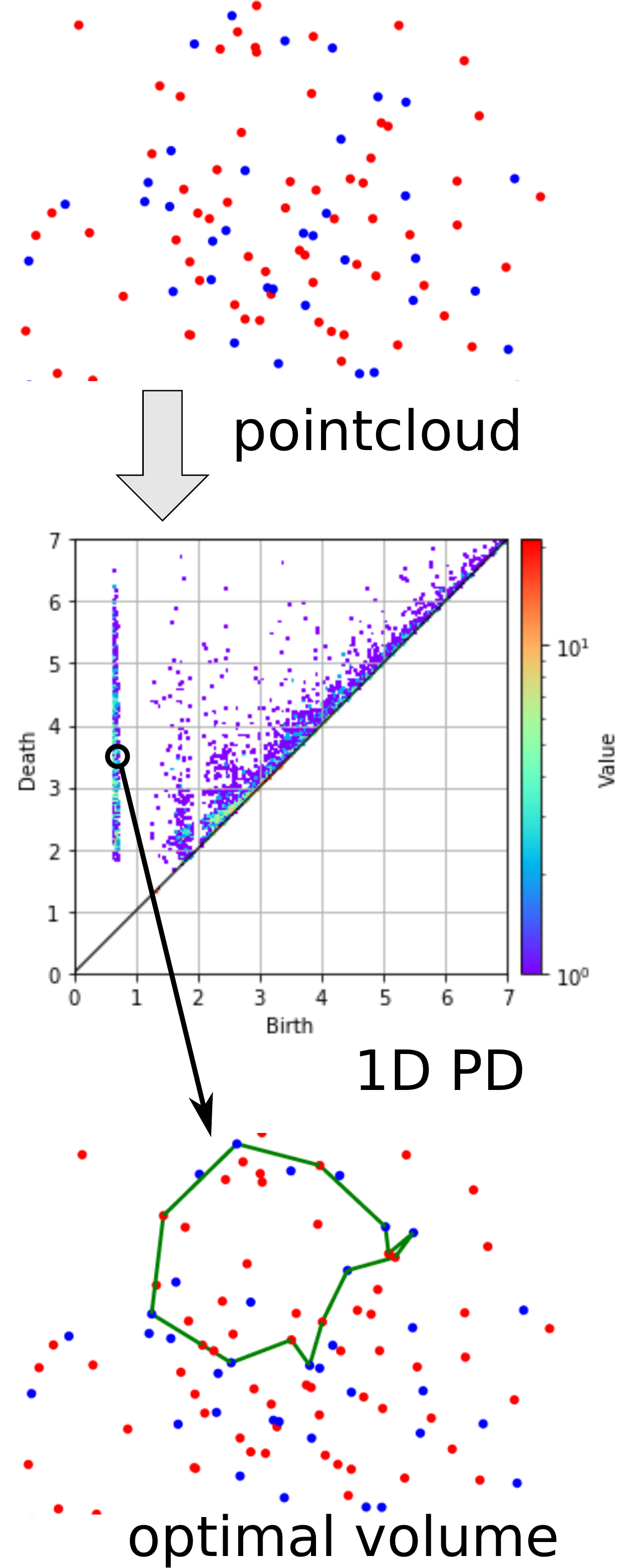}
  \caption{Inverse analysis}
  \label{fig:inverse}
\end{figure}

HomCloud has advanced inverse analysis features. Figure~\ref{fig:inverse}
shows the outline of the functionality. In a PD,
each birth-death pair corresponds to a ring or a cavity. It is very helpful
to identify the structure to analyze PDs, but the identification is not easy
since many candidates exist. Mathematical optimization is used to
select the tightest structure from the candidates\cite{Tahbaz-Salehi,oc2,eh,voc,p1cycle,stable-volume}. HomCloud already implements some methods\cite{voc,stable-volume}. We can
apply different methods depending on the purpose, input type, and performance.
The inverse analysis is available for all types of input,
but which method to use depends on the data type.

We can easily visualize the result of inverse analysis using HomCloud.
We can also output the geometric information of the result in several forms.
We can compare the output with other information for further investigation.

\subsubsection{Machine learning}

HomCloud supports machine learning using PDs.
Machine learning can find hidden patterns that are common to many PDs.
Since many machine learning methods require feature vectors or Gram matrices as input data, we need to convert PDs into vectors or matrices\cite{ph_ml_survey}.
HomCloud supports persistence image\cite{pi} for a vectorization method. 
Persistence image vectors are computed from the histograms of
diagrams. Intuitively saying,
the values of bins of the histogram are used as vector elements.
The following techniques are used to improve the performance of learning.
\begin{itemize}
\item 2D Gaussian filter to reflect the adjacency of bins
\item Weight function to reflect the importance of birth-death pairs
  depending on the distance from the diagonal
\end{itemize}

The advantage of persistence image is its simplicity. We can
intuitively understand the method. We can also convert a vector into a histogram
in reverse. The advantage is especially useful if we use linear machine learning
models such as linear regression, logistic regression, and principal component
analysis. Since the models give vectors with the same dimension as
input vectors, we can visualize the learned results in histograms.
The histograms show which birth-death pairs are important. After identifying
important birth-death pairs, we can apply inverse analysis to map
the learned results onto the original data\cite{linear-models}.



\subsubsection{Other functionalities}

HomCloud has some other utilities. For example, HomCloud supports
bottleneck and Wasserstein distances to compare PDs using HERA\cite{hera}.
You can plot the slice of a PD histogram using HomCloud.

\subsubsection{HomCloud examples}

As explained above, HomCloud has a Python interface and command-line interface.
Comprehensive tutorials (both in English and Japanese) are available
on HomCloud's website (\url{https://homcloud.dev/basic-usage.en.html}).
You can learn how to analyze pointclouds, 2D or 3D bitmap data and distance matrices in those tutorials using Jupyter notebook (\url{https://jupyter.org/}).
Tutorials for both the Python interface and command-line interface
are available. We also prepare a tutorial about PH with machine learning.

Now we introduce the Python interface. You can try the codes on jupyter notebook or other integrated development environments such as VSCode.

First, HomCloud, Numpy, and Matplotlib are imported.
\begin{lstlisting}[language=Python]
import homcloud.interface as hc
import numpy as np
import matplotlib.pyplot as plt
\end{lstlisting}

The input pointcloud data should be a NumPy 2D array. The following code
reads the data from a text file.
\begin{lstlisting}[language=Python]
pointcloud = np.loadtxt("pointcloud.txt")
\end{lstlisting}
The text file should be as follows.
These numbers are the values of the X, Y, and Z coordinates of each point.
The X, Y, and Z values are separated by spaces or tabs.
\begin{verbatim}
  5.043  -16.116   -4.787
  7.184  -16.066   -3.850
  -8.529   -1.029   -0.326
  :         :        :
\end{verbatim}
The PDs are computed as follows:
\begin{lstlisting}
hc.PDList.from_alpha_filtration(
  pointcloud,
  save_to="pointcloud.pdgm",
  save_boundary_map=True
)
\end{lstlisting}
The PDs are saved into a file ``pointcloud.pdmg''. You can
load the information of the 1D PD by the following code:
\begin{lstlisting}
pdlist = hc.PDList("pointcloud.pdgm")
pd1 = pdlist.dth_diagram(1)
\end{lstlisting}
Since ``pointcloud.pdgm'' has the information of PDs of all dimensions,
\verb|hc.PDList("pointcloud.pdgm")| returns \verb|PDList| object
and we get \verb|PD| object by \verb|pdlist.dth_diagram(1)|.

The following code will give you the list of birth times as
a NumPy array.
\begin{lstlisting}
pd1.births
# => array([0.140, 0.139, ..., 
\end{lstlisting}
Of course, \verb|deaths| is also available to get death times.
You can construct a histogram of birth-death pairs and plot it
using HomCloud.
\begin{lstlisting}
pd1.histogram(
  (-0.2, 3.0), 256
).plot(
  colorbar={"type": "log"}
)
plt.show()  
\end{lstlisting}
You can specify the range and the number of bins as above.
\verb|(-0.2, 3.0)| is the range of birth and death times,
and \verb|256| is the number of bins in both axes.
We can also modify the color bar of the plot.
Since HomCloud uses Matplotlib as a backend, you can
create good appearance figures by calling the Matplotlib functions.

You can compute an optimal volume, one important feature of HomCloud, by the following codes:
\begin{lstlisting}
pair = pd1.nearest_pair_to(0.72, 2.1)
optimal_volume = pair.optimal_volume()  
\end{lstlisting}
The first line picks up a birth-death pair nearest to $(0.72, 2.1)$ and the second line
computes the optimal volume of the pair.
You can get access to the coordinates of points on the corresponding ring or cavity as follows.
\begin{lstlisting}
optimal_volume.boundary_points()
# => [[0.24, 0.36, 0.39],
#     [0.02, 0.10, 0.68],
#       :
\end{lstlisting}
You can also get access to the edges of a ring or the faces of a cavity using
\verb|boundary| method. HomCloud also can visualize
the rings and cavities as in Fig.~\ref{fig:inverse}.

\subsubsection{Performance remarks}

Since HomCloud uses efficient software such as PHAT or Ripser as a backend, HomCloud has good performance.
You can analyze 1,000,000 points in 3D space or 300x300x300 voxel data
in five or ten minutes with a recent PC with 32GB memory.
Larger data requires a larger computer.
Some benchmark results are available at \url{https://homcloud.dev/benchmarks.html}.

\subsubsection{Future plans}

Our developer team has a plan to continuously improve HomCloud.
New features will be implemented in parallel with theoretical research.
Performance improvements and UI enhancements will be made as needed.
We will also implement new methods proposed by previous researches.
Now we have the idea to implement vectorization and kernel methods
for machine learning other than persistence image. 

\section{Concluding Remarks}\label{sec:conclusion}
PH gives a way to summarize the shape of data quantitatively and it is helpful to analyze materials data from micro-scale data to macro-scale data.
PH is now rapidly developed from theory, software, and to applications to various fields including materials science.
The collaboration between theory and applications of PH is quite active.

One example of collaboration is the study of multi-parameter persistence.
In PH, scale information is encoded in a filtration, enabling us to extract multi-scale geometric structures.
If we want to integrate another information such as noise reduction level or temporal information, it is natural to introduce another filtration axis.
Carlsson and Zomorodian\cite{multidimensional} first proposed multi-parameter persistence and \emph{rank invariant} to characterize it.
The paper also showed the theoretical difficulty of multi-parameter persistence.
Now mathematicians tackle this mathematical problem to give a better characterization\cite{Escolar2016,Kim2021,asashiba2021approximation}, and in the future, they will provide feedback for applications.

Software is essential for connecting theory to applications.
A lot of PH software has been developed according to the theoretical interests and applications of developers.
One of the software, HomCloud, is introduced in this paper.
HomCloud has some advanced features such as inverse analysis and machine learning support and has been applied to data analysis of materials.

This paper introduces PH for readers interested in the application to materials science and gives an intuitive explanation of PH.
Some review papers and textbooks\cite{comptop,ghrist,carlsson,roadmap,ehsurvay,Buchet2018,surveyalg}, including a textbook in Japanese\cite{hiraokatextbookjp}, have been written for readers who are interested in further mathematical details and algorithms.

\section*{Acknowledgments}
This work was partially sponsored by JSPS KAKENHI JP19H00834,
JP20H05884, JP18H01188, JST Presto JPMJPR1923, JST CREST Mathematics (JPMJCR15D3),  JST MIRAI Program (JPMJMI18G3), and Council for Science, Technology and Innovation (CSTI), Cross-ministerial Strategic Innovation Promotion Program (SIP) and ``Materials Integration'' for revolutionary design system of structural materials (Funding agency: JST).

\bibliography{refs} 

\begin{thebibliography}{10}

\bibitem{elz}
H.~Edelsbrunner, D.~Letscher, and A.~Zomorodian: Proceedings 41st annual
  symposium on foundations of computer science, 2000, pp. 454--463.

\bibitem{zc}
A.~Zomorodian and G.~Carlsson: Discrete {\&} Computational Geometry {\bfseries
  33} (2005) 249.

\bibitem{ceh}
D.~Cohen-Steiner, H.~Edelsbrunner, and J.~Harer: Discrete {\&} Computational
  Geometry {\bfseries 37} (2007) 103.

\bibitem{csggo}
F.~Chazal, D.~Cohen-Steiner, M.~Glisse, L.~J. Guibas, and S.~Y. Oudot:
  Proceedings of the twenty-fifth annual symposium on Computational geometry,
  2009, pp. 237--246.

\bibitem{bl}
U.~Bauer and M.~Lesnick: Proceedings of the thirtieth annual symposium on
  Computational geometry, 2014, pp. 355--364.

\bibitem{dg}
V.~De~Silva and R.~Ghrist: Algebraic \& Geometric Topology {\bfseries 7} (2007)
  339.

\bibitem{hnhemn}
Y.~Hiraoka, T.~Nakamura, A.~Hirata, E.~G. Escolar, K.~Matsue, and Y.~Nishiura:
  Proceedings of the National Academy of Sciences {\bfseries 113} (2016) 7035.

\bibitem{stfrh}
M.~Saadatfar, H.~Takeuchi, V.~Robins, N.~Francois, and Y.~Hiraoka: Nature
  communications {\bfseries 8} (2017) 1.

\bibitem{virus}
J.~M. Chan, G.~Carlsson, and R.~Rabadan: Proceedings of the National Academy of
  Sciences {\bfseries 110} (2013) 18566.

\bibitem{wei}
G.-W. Wei: Journal of Computational Physics {\bfseries 305} (2017) 276.

\bibitem{brain}
C.~Giusti, E.~Pastalkova, C.~Curto, and V.~Itskov: Proceedings of the National
  Academy of Sciences {\bfseries 112} (2015) 13455.

\bibitem{cosmology}
R.~Van De~Weygaert, G.~Vegter, H.~Edelsbrunner, B.~J. Jones, P.~Pranav,
  C.~Park, W.~A. Hellwing, B.~Eldering, N.~Kruithof, E.~P. Bos, et~al.,
  Transactions on computational science XIV, pp. 60--101. Springer, 2011.

\bibitem{eh}
E.~G. Escolar and Y.~Hiraoka: {\em Optimal Cycles for Persistent Homology Via
  Linear Programming} (Springer Japan, Tokyo, 2016), pp. 79--96.

\bibitem{linear-models}
I.~Obayashi, Y.~Hiraoka, and M.~Kimura: Journal of Applied and Computational
  Topology {\bfseries 1} (2018) 421.

\bibitem{ph_ml_survey}
C.~S. Pun, K.~Xia, and S.~X. Lee: arXiv:1811.00252 .

\bibitem{alpha}
H.~Edelsbrunner and E.~P. M\"{u}cke: ACM Trans. Graph. {\bfseries 13} (1994)
  43–72.

\bibitem{iron-ore}
M.~Kimura, I.~Obayashi, Y.~Takeichi, R.~Murao, and Y.~Hiraoka: Scientific
  reports {\bfseries 8} (2018) 1.

\bibitem{porous}
V.~Robins, M.~Saadatfar, O.~Delgado-Friedrichs, and A.~P. Sheppard: Water
  Resources Research {\bfseries 52} (2016) 315.

\bibitem{capillary-fluid}
A.~Herring, V.~Robins, and A.~Sheppard: Water Resources Research {\bfseries 55}
  (2019) 555.

\bibitem{flow-estimation}
A.~Suzuki, M.~Miyazawa, J.~M. Minto, T.~Tsuji, I.~Obayashi, Y.~Hiraoka, and
  T.~Ito: Scientific reports {\bfseries 11} (2021) 1.

\bibitem{nakamura2015persistent}
T.~Nakamura, Y.~Hiraoka, A.~Hirata, E.~G. Escolar, and Y.~Nishiura:
  Nanotechnology {\bfseries 26} (2015) 304001.

\bibitem{cheng2011atomic}
Y.~Cheng and E.~Ma: Progress in materials science {\bfseries 56} (2011) 379.

\bibitem{zallen2008physics}
R.~Zallen: {\em The physics of amorphous solids} (John Wiley \& Sons, 2008).

\bibitem{steinhardt1983bond}
P.~J. Steinhardt, D.~R. Nelson, and M.~Ronchetti: Physical Review B {\bfseries
  28} (1983) 784.

\bibitem{hansen1990theory}
J.-P. Hansen and I.~R. McDonald: {\em Theory of simple liquids} (Elsevier,
  1990).

\bibitem{delhommelle2001inadequacy}
J.~Delhommelle and P.~Milli{\'e}: Molecular Physics {\bfseries 99} (2001) 619.

\bibitem{kondic2012topology}
L.~Kondic, A.~Goullet, C.~O'Hern, M.~Kramar, K.~Mischaikow, and R.~Behringer:
  EPL (Europhysics Letters) {\bfseries 97} (2012) 54001.

\bibitem{ardanza2014topological}
S.~Ardanza-Trevijano, I.~Zuriguel, R.~Ar{\'e}valo, and D.~Maza: Physical Review
  E {\bfseries 89} (2014) 052212.

\bibitem{lim2018topology}
M.~X. Lim and R.~P. Behringer: EPL (Europhysics Letters) {\bfseries 120} (2018)
  44003.

\bibitem{hayakawa2021impact}
H.~Hayakawa et~al.: Physical Review Fluids {\bfseries 6} (2021) 033301.

\bibitem{coslovich2009dynamics}
D.~Coslovich and G.~Pastore: Journal of Physics: Condensed Matter {\bfseries
  21} (2009) 285107.

\bibitem{ormrod2021persistent}
D.~Ormrod~Morley, P.~S. Salmon, and M.~Wilson: The Journal of Chemical Physics
  {\bfseries 154} (2021) 124109.

\bibitem{kusano2017kernel}
G.~Kusano, K.~Fukumizu, and Y.~Hiraoka: The Journal of Machine Learning
  Research {\bfseries 18} (2017) 6947.

\bibitem{Onodera2020}
Y.~Onodera, S.~Kohara, P.~S. Salmon, A.~Hirata, N.~Nishiyama, S.~Kitani,
  A.~Zeidler, M.~Shiga, A.~Masuno, H.~Inoue, S.~Tahara, A.~Polidori, H.~E.
  Fischer, T.~Mori, S.~Kojima, H.~Kawaji, A.~I. Kolesnikov, M.~B. Stone, M.~G.
  Tucker, M.~T. McDonnell, A.~C. Hannon, Y.~Hiraoka, I.~Obayashi, T.~Nakamura,
  J.~Akola, Y.~Fujii, K.~Ohara, T.~Taniguchi, and O.~Sakata: NPG Asia Materials
  {\bfseries 12} (2020) 85.

\bibitem{Koyama2020}
C.~Koyama, S.~Tahara, S.~Kohara, Y.~Onodera, D.~R. Sm{\aa}br{\aa}ten, S.~M.
  Selbach, J.~Akola, T.~Ishikawa, A.~Masuno, A.~Mizuno, J.~T. Okada,
  Y.~Watanabe, Y.~Nakata, K.~Ohara, H.~Tamaru, H.~Oda, I.~Obayashi, Y.~Hiraoka,
  and O.~Sakata: NPG Asia Materials {\bfseries 12} (2020) 43.

\bibitem{hong2019medium}
S.~Hong and D.~Kim: Journal of Physics: Condensed Matter {\bfseries 31} (2019)
  455403.

\bibitem{Shimizu2021}
Y.~Shimizu, T.~Kurokawa, H.~Arai, and H.~Washizu: Scientific Reports {\bfseries
  11} (2021) 2274.

\bibitem{yoshimoto2021molecular}
Y.~Yoshimoto, S.~Sugiyama, S.~Shimada, T.~Kaneko, S.~Takagi, and I.~Kinefuchi:
  Macromolecules {\bfseries 54} (2021) 958.

\bibitem{ichinomiya2017persistent}
T.~Ichinomiya, I.~Obayashi, and Y.~Hiraoka: Physical Review E {\bfseries 95}
  (2017) 012504.

\bibitem{PhysRevResearch.2.033529}
F.~Landuzzi, T.~Nakamura, D.~Michieletto, and T.~Sakaue: Phys. Rev. Research
  {\bfseries 2} (2020) 033529.

\bibitem{hosokawa2019partial}
S.~Hosokawa, J.-F. B{\'e}rar, N.~Boudet, W.-C. Pilgrim, L.~Pusztai, S.~Hiroi,
  K.~Maruyama, S.~Kohara, H.~Kato, H.~E. Fischer, et~al.: Physical Review B
  {\bfseries 100} (2019) 054204.

\bibitem{rocks2021learning}
J.~W. Rocks, S.~A. Ridout, and A.~J. Liu: APL Materials {\bfseries 9} (2021)
  021107.

\bibitem{shirai2019microscopic}
T.~Shirai and T.~Nakamura: Journal of the Physical Society of Japan {\bfseries
  88} (2019) 074801.

\bibitem{doi:10.1021/acs.langmuir.8b01411}
V.~Lotito and T.~Zambelli: Langmuir {\bfseries 34} (2018) 7827.
\newblock PMID: 29886749.

\bibitem{carr2016energy}
J.~M. Carr, D.~Mazauric, F.~Cazals, and D.~J. Wales: The Journal of chemical
  physics {\bfseries 144} (2016) 054109.

\bibitem{Kotsugi}
T.~Yamada, Y.~Suzuki, C.~Mitsumata, K.~Ono, T.~Ueno, I.~Obayashi, Y.~Hiraoka,
  and M.~Kotsugi: Vacuum and Surface Science {\bfseries 62} (2019) 153.
\newblock [in Japanese].

\bibitem{published_papers/31169608}
Y.~Mototake, S.~Yamanaka, T.~Aoyagi, T.~Ohnishi, and K.~Fukumizu: Proceedings
  of the 2020 International Symposium on Nonlinear Theory and its Applications
  (2020).

\bibitem{twist}
C.~Chen and M.~Kerber: Proceedings 27th European Workshop on Computational
  Geometry, Vol.~11, 2011, pp. 197--200.

\bibitem{cohomology1}
V.~De~Silva, D.~Morozov, and M.~Vejdemo-Johansson: Inverse Problems {\bfseries
  27} (2011) 124003.

\bibitem{cohomology2}
V.~De~Silva, D.~Morozov, and M.~Vejdemo-Johansson: Discrete \& Computational
  Geometry {\bfseries 45} (2011) 737.

\bibitem{eirene}
G.~Henselman and R.~Ghrist: arXiv:1606.00199 .

\bibitem{multifield}
J.-D. Boissonnat and C.~Maria: Journal of Applied and Computational Topology
  {\bfseries 3} (2019) 59.

\bibitem{chunks}
U.~Bauer, M.~Kerber, and J.~Reininghaus, Topological methods in data analysis
  and visualization III, pp. 103--117. Springer, 2014.

\bibitem{dipha}
U.~Bauer, M.~Kerber, and J.~Reininghaus: 2014 proceedings of the sixteenth
  workshop on algorithm engineering and experiments (ALENEX), 2014, pp. 31--38.

\bibitem{perseus}
K.~Mischaikow and V.~Nanda: Discrete {\&} Computational Geometry {\bfseries 50}
  (2013) 330.

\bibitem{phat}
U.~Bauer, M.~Kerber, J.~Reininghaus, and H.~Wagner: Journal of Symbolic
  Computation {\bfseries 78} (2017) 76.
\newblock Algorithms and Software for Computational Topology.

\bibitem{ripser}
U.~Bauer: Journal of Applied and Computational Topology {\bfseries 5} (2021)
  391.

\bibitem{giotto-tda}
G.~Tauzin, U.~Lupo, L.~Tunstall, J.~B. P{\'e}rez, M.~Caorsi, A.~M.
  Medina-Mardones, A.~Dassatti, and K.~Hess: J. Mach. Learn. Res. {\bfseries
  22} (2021) 39.

\bibitem{jholes}
J.~Binchi, E.~Merelli, M.~Rucco, G.~Petri, and F.~Vaccarino: Electronic Notes
  in Theoretical Computer Science {\bfseries 306} (2014) 5.
\newblock Proceedings of the 5th International Workshop on Interactions between
  Computer Science and Biology (CS2Bio’14).

\bibitem{roadmap}
N.~Otter, M.~A. Porter, U.~Tillmann, P.~Grindrod, and H.~A. Harrington: EPJ
  Data Science {\bfseries 6} (2017) 17.

\bibitem{Hirata2020}
A.~Hirata, T.~Wada, I.~Obayashi, and Y.~Hiraoka: Communications Materials
  {\bfseries 1} (2020) 98.

\bibitem{YoheiONODERA201919143}
Y.~Onodera, S.~Kohara, S.~Tahara, A.~Masuno, H.~INoue, M.~Shiga, A.~Hirata,
  K.~Tsuchiya, Y.~Hiraoka, I.~Obayashi, K.~Ohara, A.~Mizuno, and O.~Sakata:
  Journal of the Ceramic Society of Japan {\bfseries 127} (2019) 853.

\bibitem{Onodera2019a}
Y.~Onodera, Y.~Takimoto, H.~Hijiya, T.~Taniguchi, S.~Urata, S.~Inaba,
  S.~Fujita, I.~Obayashi, Y.~Hiraoka, and S.~Kohara: NPG Asia Materials
  {\bfseries 11} (2019) 75.

\bibitem{PhysRevB.99.045153}
M.~Murakami, S.~Kohara, N.~Kitamura, J.~Akola, H.~Inoue, A.~Hirata, Y.~Hiraoka,
  Y.~Onodera, I.~Obayashi, J.~Kalikka, N.~Hirao, T.~Musso, A.~S. Foster,
  Y.~Idemoto, O.~Sakata, and Y.~Ohishi: Phys. Rev. B {\bfseries 99} (2019)
  045153.

\bibitem{doi:10.1098/rspa.2020.0170}
M.~Cramer~Pedersen, V.~Robins, K.~Mortensen, and J.~J.~K. Kirkensgaard:
  Proceedings of the Royal Society A: Mathematical, Physical and Engineering
  Sciences {\bfseries 476} (2020) 20200170.

\bibitem{ANDO2021142112}
I.~Ando, Y.~Mugita, K.~Hirayama, S.~Munetoh, M.~Aramaki, F.~Jiang, T.~Tsuji,
  A.~Takeuchi, M.~Uesugi, and Y.~Ozaki: Materials Science and Engineering: A
  {\bfseries 828} (2021) 142112.

\bibitem{Hong_2019}
S.~Hong and D.~Kim: Journal of Physics: Condensed Matter {\bfseries 31} (2019)
  455403.

\bibitem{minamitani2021topological}
E.~Minamitani, T.~Shiga, M.~Kashiwagi, and I.~Obayashi:  arXiv:2107.05865.

\bibitem{SUZUKI2020104550}
A.~Suzuki, M.~Miyazawa, A.~Okamoto, H.~Shimizu, I.~Obayashi, Y.~Hiraoka,
  T.~Tsuji, P.~Kang, and T.~Ito: Computers \& Geosciences {\bfseries 143}
  (2020) 104550.

\bibitem{ICHINOMIYA20202926}
T.~Ichinomiya, I.~Obayashi, and Y.~Hiraoka: Biophysical Journal {\bfseries 118}
  (2020) 2926.

\bibitem{Koseki2020}
K.~Koseki, H.~Kawasaki, T.~Atsugi, M.~Nakanishi, M.~Mizuno, E.~Naru,
  T.~Ebihara, M.~Amagai, and E.~Kawakami: npj Systems Biology and Applications
  {\bfseries 6} (2020) 40.

\bibitem{Oyama2019}
A.~Oyama, Y.~Hiraoka, I.~Obayashi, Y.~Saikawa, S.~Furui, K.~Shiraishi,
  S.~Kumagai, T.~Hayashi, and J.~Kotoku: Scientific Reports {\bfseries 9}
  (2019) 8764.

\bibitem{weighted-alpha}
H.~Edelsbrunner: Technical report, Champaign, IL, USA (1992).

\bibitem{Tahbaz-Salehi}
A.~Tahbaz-Salehi and A.~Jadbabaie: 2008 47th IEEE Conference on Decision and
  Control, Dec 2008, pp. 4170--4176.

\bibitem{oc2}
T.~K. Dey, A.~N. Hirani, and B.~Krishnamoorthy: SIAM J. Comput. {\bfseries 40}
  (2011) 1026.

\bibitem{voc}
I.~Obayashi: SIAM Journal on Applied Algebra and Geometry {\bfseries 2} (2018)
  508.

\bibitem{p1cycle}
T.~K. Dey, T.~Hou, and S.~Mandal: In R.~Marfil, M.~Calder{\'o}n,
  F.~D{\'i}az~del R{\'i}o, P.~Real, and A.~Bandera (eds), {\em Computational
  Topology in Image Context}, 2019, pp. 123--136.

\bibitem{stable-volume}
I.~Obayashi: arXiv:2109.11711 .

\bibitem{pi}
H.~Adams, T.~Emerson, M.~Kirby, R.~Neville, C.~Peterson, P.~Shipman,
  S.~Chepushtanova, E.~Hanson, F.~Motta, and L.~Ziegelmeier: Journal of Machine
  Learning Research {\bfseries 18} (2017).

\bibitem{hera}
M.~Kerber, D.~Morozov, and A.~Nigmetov: Journal of Experimental Algorithmics
  {\bfseries 22} (2017) 1.

\bibitem{multidimensional}
G.~Carlsson and A.~Zomorodian: Discrete {\&} Computational Geometry {\bfseries
  42} (2009) 71.

\bibitem{Escolar2016}
E.~G. Escolar and Y.~Hiraoka: Discrete {\&} Computational Geometry {\bfseries
  55} (2016) 100.

\bibitem{Kim2021}
W.~Kim and F.~M{\'e}moli: Journal of Applied and Computational Topology
  {\bfseries 5} (2021) 533.

\bibitem{asashiba2021approximation}
H.~Asashiba, E.~G. Escolar, K.~Nakashima, and M.~Yoshiwaki:   (2021).
\newblock arXiv:1911.01637.

\bibitem{comptop}
H.~Edelsbrunner and J.~Harer: {\em Computational topology: an introduction}
  (American Mathematical Soc., 2010).

\bibitem{ghrist}
R.~Ghrist: {\em Elementary Applied Topology} (Createspace, 2014) 1.0 ed.

\bibitem{carlsson}
G.~Carlsson: Bull. Amer. Math. Soc. {\bfseries 46} (2009) 255.

\bibitem{ehsurvay}
H.~Edelsbrunner and J.~Harer: {\em Persistent homology—a survey} (Amer. Math.
  Soc., Providence, RI, 2008), Vol. 453 of {\em Contemp. Math.}, pp. 257--282.

\bibitem{Buchet2018}
M.~Buchet, Y.~Hiraoka, and I.~Obayashi: in{\em Persistent Homology and
  Materials Informatics}, ed. I.~Tanaka (Springer Singapore, Singapore, 2018),
  pp. 75--95.

\bibitem{surveyalg}
M.~Vejdemo-Johansson: {\em Sketches of a platypus: a survey of persistent
  homology and its algebraic foundations} (Amer. Math. Soc., 2014), Vol. 620 of
  {\em Contemp. Math.}

\bibitem{hiraokatextbookjp}
Y.~Hiraoka: {\em Tanpakushitu kouzou to topology: Persistent homology gun
  nyumon (Structure of protein and topology: Introduction to persistent
  homology)} (Kyoritu Shuppan, Tokyo, 2013), [in Japanese].

\end{thebibliography}
\bibliographystyle{jpsj}

\end{document}